\documentclass[12pt]{article}

\usepackage{latexsym}
\usepackage{calc}
\usepackage{amssymb, amsmath, amsfonts, listings, mathrsfs} 

\def\argmin{\mathop{\rm argmin}}

\def\nr{\par \noindent}

\def\Def{\stackrel{\mathrm{def}}{=}}

\def\inter{{\rm int \,}}
\def\Ind{{\rm Ind}}

\def\dom{{\rm dom \,}}
\def\beq{\begin{equation}}
\def\eeq{\end{equation}}

\def\R{\mathbb{R}}
\def\E{\mathbb{E}}

\def\vp{\varphi}
\def\D{\mathscr{D}}

\def\BI{\begin{itemize}}
\def\EI{\end{itemize}}
\def\II{\item}
\newcommand{\SetEQ}{\setcounter{equation}{0}}
\newcommand{\refLE}[1]{\ensuremath{\stackrel{\eqref{#1}}{\leq}}}
\newcommand{\refEQ}[1]{\ensuremath{\stackrel{\eqref{#1}}{=}}}
\newcommand{\refGE}[1]{\ensuremath{\stackrel{\eqref{#1}}{\geq}}}

\newtheorem{theorem}{Theorem}
\newtheorem{lemma}{Lemma}
\newtheorem{corollary}{Corollary}

\newtheorem{assumption}{Assumption}
\newtheorem{definition}{Definition}

\newtheorem{example}{Example}
\newtheorem{remark}{Remark}
\newcommand{\proof}{\bf Proof: \rm \nr}
\newcommand{\qed}{\hfill $\Box$ \nr \medskip}

\def\ba{\begin{array}}
\def\ea{\end{array}}
\def\beann{\begin{eqnarray*}}
\def\eeann{\end{eqnarray*}}
\def\bea{\begin{eqnarray}}
\def\eea{\end{eqnarray}}

\def\BT{\begin{theorem}}
\def\ET{\end{theorem}}
\def\BL{\begin{lemma}}
\def\EL{\end{lemma}}
\def\BC{\begin{corollary}}
\def\EC{\end{corollary}}
\def\BE{\begin{example}}
\def\EE{\end{example}}
\def\BD{\begin{definition}}
\def\ED{\end{definition}}
\def\BR{\begin{remark}}
\def\ER{\end{remark}}
\def\BAS{\begin{assumption}}
\def\EAS{\end{assumption}}
\def\BI{\begin{itemize}}
\def\EI{\end{itemize}}

\def\BMP{\begin{minipage}{9.5cm}}
\def\EMP{\end{minipage}}
\def\MPT{\begin{minipage}{11.5cm}}
\def\EPT{\end{minipage}}

\def\la{\langle}
\def\ra{\rangle}

\def\QF{\hspace{5ex} \Box}
\def\QR{\hfill \Box}

\title{
{\normalsize CORE DISCUSSION PAPER }\\{\normalsize
2021/01}\\\vspace{10mm} \textbf{Optimization Methods for
Fully Composite Problems} \thanks{Research results presented in this paper were obtained
in the framework of ERC Advanced Grant~788368.} }

\author{Nikita Doikov \thanks{Institute of Information and Communication Technologies,
	Electronics and Applied Math. (ICTEAM), Catholic University of Louvain (UCL). E-mail:
Nikita.Doikov@uclouvain.be. ORCID: 0000-0003-1141-1625.} 
\and Yurii Nesterov
\thanks{Center for Operations Research and Econometrics (CORE),
Catholic University of Louvain (UCL). E-mail:
Yurii.Neterov@uclouvain.be. ORCID: 0000-0002-0542-8757.}}

\date{ March 20, 2021 
}

\begin{document}
\maketitle

\abstract{In this paper, 
we propose a new Fully Composite Formulation
of convex optimization problems.
It includes, as a particular case,
the problems with functional constraints,
max-type minimization problems, and problems of Composite Minimization, where the objective can have 
simple nondifferentiable components.
We treat all these formulations in a unified way,
highlighting the existence of very natural optimization schemes
of different order. We prove the global convergence rates 
for our methods under the most general conditions.
Assuming that the upper-level component of our objective function is subhomogeneous, we develop efficient modification
of the basic Fully Composite first-order and second-order Methods,
and propose their accelerated variants.
}

\vspace{10ex}\noindent

\thispagestyle{empty}

\newpage\setcounter{page}{1}

\section{Introduction}
\setcounter{equation}{0}

\vspace{1ex}\noindent
{\bf Motivation.} Development of the numerical methods for
solving different optimization problems heavily depends on
the {\em model of the problem} used by the method's
designer. In modern Optimization Theory, the diversity of
problem formulations is sufficiently big. We can speak
about problems with functional constraints, or with simple
feasible set. The problems can be posed with
differentiable or non-differentiable components. Sometimes
we speak about problems in (additive) composite form (e.g.
\cite{nesterov2013gradient}). Or, we can speak about optimization of
max-type functions (e.g. Section 2.3 in \cite{nesterov2018lectures}).

All these formulations have quite specific properties and
usually they need development of the specific methods. In
this paper, we are going to take a step back in this
picture and consider a very general problem formulation
which covers practically all variants of the existing
problem settings. The main advantage of our formulation
(we call it {\em Fully Composite Optimization Problem}) is
that for justification of the corresponding numerical
schemes we can use only very basic properties of our
objects (convexity, monotonicity). Thus, we highlight the
generic reasons for existence of the efficient methods for
many different problem classes.

As an immediate consequence of our results, we get, in particular, new
high-order methods with global linear rate of convergence
for convex minimization with functional constraints. Our
new first-, second-, and third-order methods can be
implemented in practice using the existing polynomial-time
technique \cite{nesterov2019implementable}.

\vspace{1ex}\noindent
{\bf Contents.}
In Section~\ref{sc-UCF}, we study uniformly 
convex smooth functions. We prove two new inequalities based on high-order Taylor polynomials, which provide these functions with the improved global lower bounds.
This gives us the main tool for justifying 
the global convergence rates of our methods.

In Section~\ref{sc-FComp}, we present 
our Fully Composite Optimization Framework, and give
several examples, which cover all popular composite settings. Then, in Section~\ref{sc-MBase}, we develop
basic high-order optimization methods
(starting from the first-order methods) for solving 
Fully Composite Problems.
Assuming that the smooth component
of our problem is uniformly convex of a certain degree,
we establish a global linear rate of convergence of 
the new methods.

In Section~\ref{sc-Reg}, we demonstrate that it is possible
to use a simple regularization technique within our framework.
It converts any convex problem into uniformly convex one, and
thus our basic methods can be applied to solve them.

Section~\ref{sc-Subhomo} is devoted to 
\textit{subhomogeneous} functions.
We provide the definition and list several properties of such functions.
Then we show that for subhomogeneous fully composite formulations,
the global convergence of the basic methods
holds in a more general convex setting.

We study efficient modifications of the first-order and second-order methods
for the subhomogeneous fully composite problems 
in Sections~\ref{sc-Gradient} and \ref{sc-Newton}
respectively. In particular, we establish the accelerated $O(k^{-2})$
global rate of convergence for the Fast Gradient Method~\cite{nesterov1983method},
and the same global rate for the modifications of 
the Newton's Method~\cite{nesterov2006cubic,doikov2020convex}.

In Section~\ref{sc-ContrProx}, we accelerate our fully composite 
second-order methods up to the level $O(k^{-3})$ 
using inexact contracting proximal iterations~\cite{doikov2020contracting}.

\vspace{1ex}\noindent
{\bf Notation.} In what follows, we denote by $\E$ a
finite-dimensional real vector space, and by $\E^*$ its
dual spaced composed by linear functions on $\E$. For such
a function $s \in \E^*$, we denote by $\la s, x \ra$ its
value at $x \in \E$. Using a self-adjoint
positive-definite operator $B: \E \to \E^*$ (notation $B =
B^* \succ 0$), we define the {\em
conjugate Euclidean norms}:
$$
\ba{rcl}
\| x \| & = & \la B x, x \ra^{1/2}, \quad x \in \E, \quad
\| g \|_* \; = \; \la g, B^{-1} g \ra^{1/2}, \quad g \in
\E^*.
\ea
$$

For a smooth function $f: \dom f \to \R$ with convex and
open domain $\dom f  \subseteq \E$, denote by $\nabla
f(x)$ its gradient, and by $\nabla^2 f(x)$ its Hessian
evaluated at point $x \in \dom f \subseteq \E$. Then
$$
\ba{rcl}
\nabla f(x) & \in & \E^*, \quad \nabla^2 f(x) h \; \in \;
\E^*, \quad x \in \dom f, \; h \in \E.
\ea
$$

In what follows, we often work with directional
derivatives. For $p \geq 1$, denote by
$$
\ba{c}
D^p f(x)[h_1, \dots, h_p]
\ea
$$
the directional derivative of function $f$ at $x$ along
directions $h_i \in \E$, $i = 1, \dots, p$. Note that $D^p
f(x)[ \cdot]$ is a {\em symmetric $p$-linear form}. Its
{\em norm} is defined in the standard way:
\beq\label{eq-DNorm}
\ba{rcl}
\| D^pf(x) \| & = & \max\limits_{h_1, \dots, h_p} \left\{
D^p f(x)[h_1, \dots, h_p ]: \; \| h_i \| \leq 1, \, i = 1,
\dots, p \right\}.
\ea
\eeq
For example, for any $x \in \dom f$ and $h_1, h_2 \in \E$,
we have
$$
\ba{rcl}
Df(x)[h_1] & = & \la \nabla f(x), h_1 \ra, \quad
D^2f(x)[h_1, h_2] \; = \; \la \nabla^2 f(x) h_1, h_2 \ra.
\ea
$$
Thus, for the Hessian, our definition corresponds to the
{\em spectral norm} of self-adjoint  linear operator
(maximal module of all eigenvalues computed with respect
to operator $B$).

If all directions $h_1, \dots, h_p$ are the same, we apply
the notation $D^p f(x)[h]^p$, $h \in \E$.
Then, Taylor approximation of function $f(\cdot)$ at $x
\in \dom f$ can be written as follows:
\beq\label{def-Taylor}
\ba{rcl}
f(y) & = & \Omega_p(f, x; y)  + o(\|y - x\|^p), \quad y
\in \dom
f,\\
\\
\Omega_{p}(f, x; y) & \Def & f(x) + \sum\limits_{k=1}^p {1
\over k!} D^k f(x)[y-x]^k, \quad y \in \E.
\ea
\eeq
Note that in general, we have (see, for example,
Appendix~1 in \cite{nesterov1994interior})
\beq\label{eq-DNorm1}
\ba{rcl}
\| D^pf(x) \| & = & \max\limits_{h} \left\{ \Big| D^p
f(x)[h]^p \Big|: \; \| h \| \leq 1 \right\}.
\ea
\eeq
Similarly, since for $x, y \in \dom f$ being fixed, the
form $D^pf(x)[\cdot, \dots, \cdot] - D^pf(y)[\cdot, \dots,
\cdot]$ is $p$-linear and symmetric, we also have
\beq\label{eq-DNorm2}
\ba{rcl}
\| D^pf(x) - D^pf(y) \| & = & \max\limits_{h} \left\{
\Big| D^p f(x)[h]^p - D^pf(y)[h]^p\Big|: \; \| h \| \leq 1
\right\}.
\ea
\eeq

In this paper, we consider functions from the problem
classes ${\cal F}_p$, which are convex and $p$ times continuously
differentiable on $\E$. Denote by $L_p$ the uniform bound
for the Lipschitz constant of $p$th derivative:
\beq\label{def-LP}
\ba{rcl}
\| D^p f(x) - D^pf(y) \| & \leq & L_p \| x - y \|, \quad
x, y \in \dom f,
\quad p \geq 1.
\ea
\eeq
Sometimes, if an ambiguity could arise, we use notation
$L_p(f)$.

Assuming that $f \in {\cal F}_{p}$ and $L_{p} < +\infty$,
by the standard integration arguments we can bound the
residual between function value and its Taylor
approximation:
\beq\label{eq-BoundF}
\ba{rcl}
| f(y) - \Omega_{p}(f, x; y) | & \leq & {L_{p} \over (p+1)!}
\| y - x \|^{p+1}, \quad x, y \in \dom f.
\ea
\eeq

\section{Uniform Convexity of Smooth
Functions}\label{sc-UCF}
\SetEQ

Let us couple our smoothness assumption with {\em
uniform convexity} of certain degree. 
Namely, let us assume
that for $p \geq 1$ there exists a constant $\sigma_{p+1}(f) > 0$ such
that
\beq\label{def-UC}
\ba{rcl}
\la \nabla f(y) - \nabla f(x) , y - x \ra & \geq &
\sigma_{p+1}(f) \| y - x \|^{p+1}, \quad x, y \in \dom f.
\ea
\eeq
By simple integration, this inequality ensures the
following functional growth:
\beq\label{eq-FGrow}
\ba{rcl}
f(y) & \geq & f(x) + \la \nabla f(x) , y - x \ra +
{\sigma_{p+1}(f) \over p+1} \| y - x \|^{p+1}, \quad x, y
\in \dom f.
\ea
\eeq

Let us consider the uniformly convex
functions of degree $p + 1$, whose
$p$th derivative is Lipschitz continuous.
We introduce the following constant:
$$
\ba{rcl}
\gamma_p(f) & \Def & \frac{\sigma_{p + 1}(f)}{L_p(f)},
\ea
$$
called the {\em condition number of degree
	$p$} of function $f$. Combining (\ref{eq-BoundF}) and (\ref{eq-FGrow}), we
get
$$
\ba{rcl}
\frac{\sigma_{p+1}(f)}{p + 1} \| y - x \|^{p+1} & \leq &
\sum\limits_{k=2}^p {1 \over k!} D^k f(x)[y-x]^k + {L_{p}
\over (p+1)!} \| y - x \|^{p+1}, \quad x, y \in \dom f.
\ea
$$
Thus, in the case of unbounded domain, we have
\beq\label{eq-GBound}
\ba{rcl}
\gamma_p(f)  & \leq & {1 \over p!}.
\ea
\eeq

Let us prove now the main inequalities of our problem class.
For $\alpha \geq 0$, denote
\beq\label{def-Beta}
\ba{rcl}
\beta_p(f,\alpha) & = & {\left( p! \, \gamma_p(f)
\right)^{1 \over p} \over (1+\alpha)^{1 \over p} +
\left(p! \, \gamma_p(f)\right)^{1 \over p}} \,.
\ea
\eeq
\BT\label{th-Main}
For any $p \geq 1$, $\alpha \geq 0$, and all $x, y \in \dom f$, we have
\beq\label{eq-GUGrow}
\ba{rcl}
\la \nabla f(y) - \nabla f(x) , y - x \ra & \geq &
\sum\limits_{k=2}^p {\beta^{k-1} \over (k-1)!} D^k
f(x)[y-x]^k + {  \alpha L_p \beta^{p} \over p!} \| y - x
\|^{p+1},
\ea
\eeq
\beq\label{eq-FUGrow}
\ba{rcl}
f(y) & \geq & f(x) + \la \nabla f(x) , y - x \ra +
\sum\limits_{k=2}^p {\beta^{k-1}\over k!} D^k f(x)[y-x]^k
+ { \alpha L_p  \beta^{p} \over (p+1)!}
 \| y - x \|^{p+1},
\ea
\eeq
where $0 \leq \beta \leq \beta_p(f,\alpha)$.
\ET
\proof
Let us fix some $x, y \in \dom f$ and consider $z_t \Def x +
t(y-x)$, $0 \leq t \leq 1$. Then,
$$
\ba{rcl}
\la \nabla f(y), y - x \ra & = & {1 \over 1 - t } \la
\nabla f(y), y - z_t \ra \\
\\
& \refGE{def-UC} & {1 \over 1 - t } \la \nabla f(z_t), y
- z_t \ra + (1-t)^p \sigma_{p+1}(f) \| y - x \|^{p+1}\\
\\
& = & \la \nabla f(z_t), y - x \ra + (1-t)^p
\sigma_{p+1}(f) \| y - x \|^{p+1}.
\ea
$$
Consider now the function $\phi(t) = \la \nabla f(z_t), y
- x \ra$. Then, by Taylor's formula, we have
$$
\ba{rcl}
\phi(t) & = & \phi(0) + \sum\limits_{k=1}^{p-1} {t^k \over
k!} \phi^{(k)}(0) + {1 \over (p-1)!} \int\limits_0^t
(t-\lambda)^{p-1} \phi^{(p)}(\lambda) d \lambda\\
\\
& = & \sum\limits_{k=1}^{p} {t^{k-1} \over (k-1)!} D^k
f(x)[y-x]^k + {1 \over (p-1)!} \int\limits_0^t
(t-\lambda)^{p-1} D^{p+1}f(z_{\lambda})[y-x]^{p+1} d
\lambda\\
\\
& \geq & \sum\limits_{k=1}^{p} {t^{k-1} \over (k-1)!} D^k
f(x)[y-x]^k - {t^p \over p!} L_p \| y - x \|^{p+1}.
\ea
$$
Adding these two inequalities, we get
$$
\ba{rcl}
\la \nabla f(y) - \nabla f(x), y - x \ra & \geq &
\sum\limits_{k=2}^{p} {t^{k-1} \over (k-1)!} D^k
f(x)[y-x]^k \\
\\
& & + \left( (1-t)^p \sigma_{p+1}(f) - {t^p \over p!}
L_p\right) \| y - x \|^{p+1}.
\ea
$$
Let us choose now $t$ from the inequality
$$
\ba{rcl}
(1-t)^p \sigma_{p+1}(f) - {1 \over p!} t^p L_p & \geq &
{\alpha \over p!} t^p L_p \quad \Leftrightarrow
\quad {p! \over 1+\alpha} \gamma_p(f) \; \geq \; \left({ t \over 1 - t} \right)^p.
\ea
$$
Then it is enough to take $t \overset{\eqref{def-Beta}}{\leq}
\beta_p(f,\alpha)$. Hence, inequality (\ref{eq-GUGrow}) is
proved.

The remaining inequality (\ref{eq-FUGrow}) can be proved
by integration. Indeed
$$
\ba{rl}
& f(y) - f(x) - \la \nabla f(x), y - x \ra \\
\\
= & \int\limits_0^1 {1 \over \tau} \la \nabla f(x +
\tau(y-x) - f(x), \tau(y-x) \ra d \tau\\
\\
\\\overset{\eqref{eq-GUGrow}}{\geq} & \int\limits_0^1
\left(\sum\limits_{k=2}^p {\beta^{k-1} \tau^{k-1} \over (k-1)!}
D^k f(x)[y-x]^k + { \alpha L_p(f) \beta^{p} \tau^p \over
p!}
\| y - x \|^{p+1} \right) d \tau\\
\\
= & \sum\limits_{k=2}^p {\beta^{k-1} \over k!} D^k
f(x)[y-x]^k + { \alpha L_p(f) \beta^{p}  \over (p+1)!} \|
y - x \|^{p+1}. \QR
\ea
$$

\BR\label{rm-Conv}
For $\alpha \geq p$, the right-hand side of inequality
(\ref{eq-FUGrow}) is convex in $y$. Indeed, let us
introduce new variables $z = x + \beta(y-x)$. Then this
right-hand side is transformed to the following function:
$$
\ba{c}
f(x) + {1 \over \beta} \left[ \la \nabla f(x) , z - x \ra
+ \sum\limits_{k=2}^p {1\over k!} D^k f(x)[z-x]^k + {
\alpha L_p  \over (p+1)!}
 \| z - x \|^{p+1} \right].
\ea
$$
Since $\alpha \geq p$, it is convex in $z$ in view of
Theorem 1 in \cite{nesterov2019implementable}.
\ER

For the optimization schemes developed in this paper,
inequality (\ref{eq-FUGrow}) serves as the main
justification tool. Let us present now the general model
of our optimization problems.

\section{Fully Composite Optimization
Problem}\label{sc-FComp}
\SetEQ

Let $F(\cdot, \cdot)$ be a function from $
\E \times \R^m$ to $\R \cup \{+\infty\}$. 
Hence, $\dom F = \{ (x, u) \in \E \times \R^m \; : \; F(x, u) < +\infty  \}$.
For each $x \in \E$, denote
$$
\ba{rcl}
{\cal D}(x) & = & \{ u \in \R^m: \; (x,u) \in \dom F\}.
\ea
$$
Our assumptions on function $F$ are as follows.
\BAS\label{ass-F}
Function $F(\cdot,\cdot)$ is closed and convex on its
domain. Moreover, for any $x \in \E$ with ${\cal D}(x)
\neq \emptyset$, function $F(x,u)$ is closed, convex
and monotone in $u \in {\cal D}(x)$.
\EAS

Consider now a vector function $f(x) = (f_1(x), \dots,
f_m(x))^T: \dom f \to \R^m$.
\BAS\label{ass-ff}
All components of function $f$ are closed and convex.
\EAS
In our framework, all information about function $f$ can
be collected by the calls of oracle of certain degree.

Let us call {\em fully composite} the following
optimization problem:
\beq\label{prob-Comp}
\fbox{$\; \vp^* \; = \; \min\limits_{x \in \dom \vp}
\left\{ \; \vp(x) \Def F(x,f(x)) \; \right\}, \;$}
\eeq
where $\dom \vp = \{ x \in \dom f \; : \; (x,f(x)) \in \dom F
\}$. We denote by $x^{*}$ a solution to problem~\eqref{prob-Comp}:
$\vp(x^{*}) = \vp^{*}$, assuming that it exists.

Of course, problem (\ref{prob-Comp}) is tractable only if
the function $F$ is simple. This is our third assumption.
\BAS\label{ass-Simp}
Structure of function $F$ is simple enough for allowing an
efficient solution of some auxiliary optimization problems
based on approximations of function $f$.
\EAS

We will see soon what kind of auxiliary problems with
function $F$ we need to solve. At this moment, let us give
several examples of fully composite optimization problems.
\begin{enumerate}
\II
{\em Optimization with functional constraints}. Consider
the following problem:
\beq\label{prob-Cons}
\min\limits_{x \in Q} \{ f_1(x): \; f_i(x) \leq 0, \; i =
2, \dots, m \},
\eeq
where $Q$ is a closed convex set and $f$ satisfies
Assumption \ref{ass-ff}. Then this problem can be written in
form (\ref{prob-Comp}) with
$$
\ba{rcl}
F(x,u) & = & u^{(1)} + \sum\limits_{i=2}^m \Ind_{\R_-}
(u^{(i)}) + \Ind_Q(x),
\ea
$$
where $\Ind_X(\cdot): \E \to \{0, +\infty\}$ is the
indicator function of the set $X \subseteq \E$.
\II
{\em Additive composite minimization} \cite{nesterov2013gradient}.
Consider the following minimization problem:
\beq\label{prob-AComp}
\min\limits_x \{ f_1(x) + \psi(x) \},
\eeq
where $f(x) \equiv \{ f_1(x) \}$ satisfies Assumption
\ref{ass-ff} and $\psi(\cdot)$ is a simple closed convex
function. Then we can take
$$
\ba{rcl}
F(x,u) & = & u^{(1)} + \psi(x).
\ea
$$
\II
{\em Functional composite minimization} (e.g.~\cite{nesterov1983method,nesterov1989effective}). 
Minimization problem
\beq\label{prob-FunComp}
\min\limits_x F(f(x)),
\eeq
where $F$ is a closed convex monotone function with $\dom
F \subseteq \R^m$, and $f$ satisfies Assumption
\ref{ass-ff}, is clearly in the form (\ref{prob-Comp}).
\end{enumerate}

Note that the Taylor polynomial (\ref{def-Taylor}) is
defined in terms of directional derivatives. Therefore we
can extend its meaning onto the vector functions without
changing notation. Similarly, we will use the following
constant vectors
$$
\ba{rcl}
L_p(f) & = & (L_p(f_1), \dots, L_p(f_m))^T,\\
\\
\sigma_{p+1}(f) & = & (\sigma_{p+1}(f_1), \dots,
\sigma_{p+1}(f_m))^T,\\
\\
\gamma_p(f) & = & (\gamma_p(f_1), \dots, \gamma_p(f_m))^T.
\ea
$$
Denote $\hat \beta_p(f) = \min\limits_{1 \leq i \leq m}
\beta_p(f_i,p)$. Then inequality (\ref{eq-FUGrow}) can be
rewritten in a vector form:
\beq\label{eq-FVGrow}
\ba{rcl}
f(y) & \geq & f(x) + \sum\limits_{k=1}^p {\beta^{k-1}\over
k!} D^k f(x)[y-x]^k + { p L_p(f) \beta^{p} \over (p+1)!}
 \| y - x \|^{p+1},
\ea
\eeq
for all $\beta \in [0, \hat \beta_p(f))$. The right-hand
side of this inequality provides us with the auxiliary
problem we need to solve at each iteration of our schemes:
\beq\label{eq-Aux}
\ba{rcl}
F\left(y, f(x) + \sum\limits_{k=1}^p {\beta^{k-1}\over k!}
D^k f(x)[y-x]^k + { p L_p(f) \beta^{p} \over (p+1)!} \| y
- x \|^{p+1} \right) & \to & \min\limits_y,
\ea
\eeq
where $\beta = \hat \beta_p(f)$. We explain the sense of
this operation in the next section.

\section{Basic High-order Optimization
Methods}\label{sc-MBase}
\SetEQ

Let $\bar x \in \dom \vp$. For inequality
(\ref{eq-FUGrow}), let us choose $\alpha = p$ and $\beta
\in (0,\hat \beta_p(f)]$. Then
$$
\ba{rcl}
 & & \!\!\!\!\!\!\!\!\!
 \!\!\!\!\!\!\!\!\!\!\!
 (1 - \beta) \vp(\bar x) + \beta \vp^* 
 \;\; = \;\;
\min\limits_{v \in \dom \vp} \Big[ (1-\beta)F(\bar x,
f(\bar x)) + \beta
F(v, f(v)) \Big]\\
\\
& \geq &  \min\limits_{v \in \dom \vp} F\Big( (1-\beta)\bar x
+ \beta v,
(1-\beta) f(\bar x) + \beta f(v) \Big)\\
\\
& \refGE{eq-FVGrow} & \min\limits_{v \in \dom \vp} F\Bigl(
(1-\beta)\bar x + \beta v, f(\bar x) + \sum\limits_{k=1}^p
{\beta^k\over k!} D^k f(\bar x)[v-\bar x]^k + { p L_p(f)
\beta^{p+1} \over (p+1)!} \| v - \bar x \|^{p+1} \Bigr)\\
\\
& = & \min\limits_{{y = \bar x + \beta(v-\bar x)} \atop {v \in
\dom \vp}} F\Bigl( y, \; \Omega_p(f, \bar{x}; y) + { p L_p(f) \over
(p+1)!} \| y - \bar x \|^{p+1} \Bigr) 
\;\; \Def \;\; \tilde
M^*_{p,\beta}(\bar x).
\ea
$$
By Remark \ref{rm-Conv}, the second argument in the
objective function of the latter problem is a
component-wise convex function. Hence, by
Assumption \ref{ass-F}, this objective is convex in $y$.

Let us look at the solution of the above minimization
problem, that is
$$
\ba{rcl}
\tilde y^*_{p,\beta}(\bar x) & \Def & \argmin\limits_y
\Big\{ F\left( y, \; \Omega_p(f, \bar{x}; y) + { p L_p(f) \over
(p+1)!}
\| y - \bar x \|^{p+1} \right) :\; \bar x + {1 \over \beta}(y - \bar x) \in \dom
\vp \Big\}.
\ea
$$
Note that in view of Assumption \ref{ass-F}, we have
$$
\ba{rcl}
\tilde M^*_{p,\beta}(\bar x) & = & F\left( \tilde
y^*_{p,\beta}(\bar x), \;
\Omega_p(f, \bar{x}; \tilde y^*_{p, \beta}(\bar x) ) + { p L_p(f) \over (p+1)!} \| \tilde
y^*_{p,\beta}(\bar x) - \bar x \|^{p+1} \right)\\
\\
& \refGE{eq-BoundF} &  F\left(\tilde y^*_{p,\beta}(\bar x),
f(\tilde y^*_{p,\beta}(\bar x) ) \right) \; = \;
\vp\left(\tilde y^*_{p,\beta}(\bar x)\right).
\ea
$$

Thus, we can estimate now the rate of convergence of the
following method.
\beq\label{met-MRest}
\ba{|l|}
\hline\\
\quad \mbox{\bf Restricted $p$th order Basic Method} \quad\\
\\
\hline\\
\ba{l}
\mbox{{\bf Choose} $x_0 \in \dom \vp$ and $\beta \in (0, \hat \beta_p(f)]$.}\\
\\
\mbox{\bf For $k \geq 0$ iterate:} \quad x_{k+1} = \tilde
y^*_{p,\beta}(x_k).\\
\\
\ea\\
\hline
\ea
\eeq

We have proved for this method the following theorem.
\BT\label{th-Rest}
Let sequence $\{ x_k \}_{k \geq 0}$ be generated by the
method (\ref{met-MRest}). Then for all $k \geq 0$ we have
\beq\label{eq-RateR}
\ba{rcl}
\vp(x_k) - \vp^* & \leq & (1-\beta)^k (\vp(x_0) - \vp^*).
\ea
\eeq
\ET
Therefore, the rate of convergence is linear,
and the contraction parameter $\beta$ can reach the condition number.

Note that method (\ref{met-MRest}) could move with bigger
steps. Indeed,
$$
\ba{rcl}
\tilde M^*_{p,\beta}(\bar x) & \geq &  \min\limits_{y \in
\dom \vp} F\left( y, \; \Omega_p(f, \bar{x}; y) + { p L_p(f) \over
(p+1)!} \| y - \bar x \|^{p+1} \right) 
\;\; \Def \;\; M^*_{p}(\bar x).
\ea
$$
Denote
$$
\ba{rcl}
y^*_{p}(\bar x) & \Def & \argmin\limits_{y \in \dom
\vp} F\left( y, \; \Omega_p(f, \bar{x}; y) + { p L_p(f) \over (p+1)!}
\| y - \bar x \|^{p+1} \right).
\ea
$$
By the same reasons as before, $M^*_{p}(\bar x) \geq
\vp\left(y^*_{p}(\bar x)\right)$. Hence, we can
estimate the rate of convergence of the following method.
\beq\label{met-MFull}
\ba{|l|}
\hline\\
\quad \mbox{\bf Full-step $p$th order Basic Method} \quad\\
\\
\hline\\
\ba{l}
\mbox{{\bf Choose} $x_0 \in \dom \vp$.}\\
\\
\mbox{\bf For $k \geq 0$ iterate:} \quad x_{k+1} =
y^*_{p}(x_k).\\
\\
\ea\\
\hline
\ea
\eeq
\BT\label{th-Full}
Let sequence $\{ x_k \}_{k \geq 0}$ be generated by
(\ref{met-MFull}). Then for all $k \geq 0$ we have
\beq\label{eq-RateF}
\ba{rcl}
\vp(x_k) - \vp^* & \leq & (1- {\hat \beta}_p(f) )^k (\vp(x_0) - \vp^*)
\ea
\eeq
\ET

In view of potentially bigger steps, method
(\ref{met-MFull}) is often faster in practice. 

\BE
Let us look at 
implementation of method
(\ref{met-MFull}) for a particular class of
optimization problems with functional 
constraints~\eqref{prob-Cons} with $Q = \E$. This is
$$
\ba{rcl}
\min\limits_{x \in \E} \{ f_1(x) : \; 
f_i(x) \leq 0, \; i = 2, \dots, m \}.
\ea
$$
Then, for $p = 1$, each iteration of the method~\eqref{met-MFull}
can be represented as follows:
$$
\ba{rcl}
y_1^{*}(\bar{x}) & = & \bar{x} \; - \; 
\biggl( \; \sum\limits_{i = 1}^m \lambda^{(i)}_{*} L_1(f_i) B \biggr)^{-1}
g(\lambda_{*}),
\ea
$$
where 
$g(\lambda) \Def \sum\limits_{i = 1}^m \lambda^{(i)} \nabla f_i(\bar{x})$,
and
$\lambda_{*} \in \R^m_+$ is a solution to the corresponding \textit{dual} problem 
\beq \label{Dual1}
\ba{c}
\max\limits_{\substack{\lambda \in \R^m_{+} \\ \lambda^{(1)} = 1}}
\biggl\{ \;
\sum\limits_{i = 1}^m \lambda^{(i)} f_i(\bar{x})
- \frac{1}{2\sum_{i = 1}^m \lambda^{(i)} L_1(f_i)} 
\| g(\lambda) \|_{*}^2
\; \biggr\}.
\ea
\eeq

On the other hand, for $p = 2$, one iteration of the method~\eqref{met-MFull} 
is as follows:
$$
\ba{rcl}
y_2^{*}(\bar{x})
& = &
\bar{x} - H(\lambda_{*}, \tau_{*})^{-1}
g(\lambda_{*}),
\ea
$$
with operator
$H(\lambda, \tau) \Def
\sum\limits_{i = 1}^m \lambda^{(i)} \nabla^2 f(\bar{x}) + \tau B$.
The optimal $\lambda_{*} \in \R^m_+$ and $\tau_{*} \in \R_+$ can be computed
from the following concave optimization problem
\beq \label{Dual2}
\ba{c}
\max\limits_{\substack{\lambda \in \R_{+}^m, \tau \in \R_{+} 
			 \\ \lambda^{(1)} = 1}}
\biggl\{ \;
\sum\limits_{i = 1}^m \lambda^{(i)} f_i(\bar{x})
- \frac{\tau^3}{6 [ \sum_{i = 1}^m \lambda^{(i)} L_2(f_i) ]^2}
- \frac{1}{2}\la H(\lambda, \tau)^{-1} g(\lambda), g(\lambda) \ra
\; \biggr\}.
\ea
\eeq
Note that typically the dimension of the problems~\eqref{Dual1} and \eqref{Dual2}
is not big. Hence, they can be solved, for example, by the
Interior-Point Methods~\cite{nesterov1994interior} very efficiently.
\EE

\section{General Regularization Scheme}\label{sc-Reg}
\SetEQ

In the previous sections, we discussed two methods for
solving problem (\ref{prob-Comp}) under assumption of
uniform convexity of functional components: $\hat
\beta_p(f) > 0$. If this assumption is not valid, we still
can apply methods of Section~\ref{sc-MBase} to a special
{\em regularized} problem.

Let us present a general regularization framework for
fully composite problem (\ref{prob-Comp}). For that, we
use a component-wise convex regularizing vector function
$$
\ba{rcl}
d(x) & = & (d_1(x), \dots, d_m(x))^T: \; \dom f \to \R^m.
\ea
$$
It is related to the starting point of our process $x_0
\in \dom \vp$ in the following way:
\beq\label{eq-DefDa}
\ba{rcl}
d(x_0) & = & f(x_0),
\ea
\eeq
\beq\label{eq-DefDb}
\ba{rcl}
d(x) \; \geq \; f(x), \quad \forall x \in \dom f.
\ea
\eeq
We discuss a simple possibility for choosing such a regularizer in the end of this section.

Let $\mu \geq 0$ be a regularizing parameter. Define the following regularized function:
\beq\label{eq-DefFM}
\ba{rcl}
\vp_{\mu}(x) & = & F(x, (1-\mu)f(x) + \mu d(x)), \quad x
\in \dom \vp.
\ea
\eeq
It is convenient to assume that the
function $d(\cdot)$ satisfies the following
assumption.
\BAS\label{ass-D}
Vector function $d(x) - f(x)$ is component-wise convex on
$\dom f$.
\EAS
In this case, the regularized function $\vp_{\mu}(\cdot)$
is convex for any $\mu \geq 0$.

Note that $\vp_{\mu}(x_0) \refEQ{eq-DefDa} \vp(x_0)$.
Clearly, for all $x \in \dom \vp$ we have
\beq\label{eq-FUp}
\ba{rcl}
\vp(x) & \overset{\eqref{eq-DefDb}}{\leq} & \vp_{\mu}(x).
\ea
\eeq

Our regularized problem looks now as follows:
\beq\label{prob-CompM}
\fbox{$\; \vp^*_{\mu} \; = \; \min\limits_{x \in \dom \vp}
	\; \vp_{\mu}(x).\;$}
\eeq
At this moment, let us assume that we are able to generate
an approximate solution to this perturbed problem by one
of the methods of Section~\ref{sc-MBase}. Namely, assume
that, for certain $\delta > 0$, we have a point $\bar x \in \dom \vp$, satisfying the
following inequality:
\beq\label{eq-AppF}
\ba{rcl}
\vp_{\mu}(\bar x) - \vp^*_{\mu} & \leq & \delta \left(
\vp_{\mu}(x_0) - \vp^*_{\mu} \right) \; \refEQ{eq-DefDa}
\; \delta \left( \vp(x_0) - \vp^*_{\mu} \right) \\
\\
& \refLE{eq-FUp} & \delta \left( \vp(x_0) - \vp^* \right).
\ea
\eeq
We need to understand now how
good is this point for our initial problem
(\ref{prob-Comp}).

In order to answer this question, we need to introduce a
{\em local measure} for non-negative vectors $g \in
\R^m_+$ with respect to some point $x \in \dom \vp$ and a
functional level $A$:
\beq\label{def-Xi}
\ba{rcl}
\xi_A(x; g) & = & \min\limits_{\lambda > 0} \left\{
\lambda: \; f(x) + {1 \over \lambda} g \in {\cal D}(x), \;
F\Big(x, f(x) + {1 \over \lambda} g\Big) \leq A \right\}.
\ea
\eeq
Clearly, this measure is well defined at least at all
points $x \in \dom \vp$ with $\vp(x) < A$.

\BL\label{lm-Contr}
Let for some points $x_0$ and $\bar x$ from $\dom \vp$ we
have $\vp(x_0) \leq A$ and $\vp (\bar x) \leq A$. Then,
for any $g \in \R^m_+$ and any coefficient $\tau \in
[0,1]$ we have
\beq\label{eq-Contr}
\ba{rcl}
\xi_A((1-\tau) x_0 + \tau\bar x; g) & \leq & {1 \over
	1-\tau}\; \xi_A(x_0; g).
\ea
\eeq
\EL
\proof
Let $f(x_0) + {1 \over \lambda}
g \in {\cal D}(x_0)$ for some $\lambda > 0$, and $F(x_0, f(x_0) + {1 \over
	\lambda} g) \leq A$. Since $f(\bar x) \in {\cal D}(\bar
x)$, we have
$$
\ba{rcl}
(1-\tau)(f(x_0) + {1 \over \lambda}g) + \tau f(\bar x) &
\in & {\cal D}((1-\tau)x_0 +\tau \bar x).
\ea
$$
By convexity, $(1-\tau)f(x_0) + \tau f(\bar x) \geq f((1-\tau)x_0
+ \tau \bar x)$, and we conclude that
$$
\ba{rcl}
f((1-\tau)x_0 + \tau \bar x)) + {1 - \tau \over \lambda}g
& \in & {\cal D}((1-\tau)x_0 + \tau \bar x).
\ea
$$
At the same time,
$$
\ba{rl}
& F((1-\tau)x_0 +\tau \bar x, f((1-\tau)x_0 + \tau \bar x))
+ {1 - \tau \over \lambda}g)\\
\\
\leq & F((1-\tau)x_0 +\tau \bar x, (1-\tau)( f(x_0) + {1
	\over \lambda}g)+ \tau f(\bar x)) )\\
\\
\leq & (1-\tau)F(x_0, f(x_0) + {1 \over \lambda}g) + \tau
F(\bar x, f(\bar x)) \; \leq \; A.
\ea
$$
Thus, $\xi_A((1-\tau) x_0 + \tau \bar x; g) \leq {\lambda
	\over 1 - \tau}$.
\qed

Denote $g_* = d(x^*) - f(x^*)$ and $\xi^*_0 = \xi_A(x_0;
g_*)$.
\BL\label{lm-Reg}
Let $A \geq \vp(x_0)$ and the regularizing function
$d(\cdot)$ satisfy Assumption \ref{ass-D} and conditions
(\ref{eq-DefDa}), (\ref{eq-DefDb}). Assume that 
parameters $\alpha$ and $\tau$ from $[0,1]$ are chosen as follows:
\beq\label{eq-Par}
\ba{rcl}
{\tau \over 1 - \tau} & \leq & {\alpha \over \mu
	\,\xi^*_0}.
\ea
\eeq
Then
\beq\label{eq-Reg}
\ba{rcl}
\vp^*_{\mu} - \vp^* & \leq & (1-\tau)(1-\alpha)(\vp(x_0) -
\vp^*) + \alpha (A - \vp^*).
\ea
\eeq
\EL
\proof
Let us fix a point $x \in \dom \vp$ with $\vp(x) \leq A$.
Denote by $g_x = d(x) - f(x) \geq 0$. Then, for any
$\alpha \in [0,1]$ we have
$$
\ba{c}
\vp_{\mu}^* \; \leq \; \vp_{\mu}(x) \; = \; F(x, f(x) +
\mu g_x) \; = \; F\left(x, (1-\alpha) f(x) + \alpha \left(
f(x) +
{\mu \over \alpha} g_x \right) \right)\\
\\
\leq \; (1-\alpha) \vp(x) + \alpha F\left(x, f(x) + {\mu
	\over \alpha} g_x \right).
\ea
$$
Let us choose now $x = (1-\tau) x_0 + \tau x^*$ with
arbitrary $\tau \in [0,1]$. Then in view of Assumption
\ref{ass-D}, we have
\beq\label{eq-GPlus}
\ba{rcl}
g_x & \leq & (1-\tau) (d(x_0) - f(x_0)) + \tau (d(x^*)-
f(x^*)) \\
\\
& \refEQ{eq-DefDa} & \tau (d(x^*)- f(x^*)) \; = \; \tau
g_*.
\ea
\eeq
Hence,
$\vp_{\mu}^* \leq (1-\tau)(1-\alpha)\vp(x_0) + \tau
(1-\alpha)\vp^* + \alpha F\left(x, f(x) + {\mu \tau \over
	\alpha} g_* \right)$.
Let our parameters satisfy inequality ${\mu
	\tau  \over \alpha(1-\tau)} \leq {1 \over \xi^*_0}$. Then
$$
\ba{rcl}
{\mu \tau \over \alpha} & \leq & {1 - \tau \over \xi^*_0}
\; \refLE{eq-Contr} \; {1 \over \xi_A(x; g^*)}.
\ea
$$
This means that $F\left(x, f(x) + {\mu \tau \over \alpha}
g_* \right) \leq A$, and we get
$$
\ba{rcl}
\vp^*_{\mu} - \vp^* & \leq & (1-\tau)(1-\alpha)(\vp(x_0) -
\vp^*) + \alpha (A - \vp^*). \QF
\ea
$$

Let us put in (\ref{eq-Reg}) the
best values of parameters. If $\tau$ satisfies
(\ref{eq-Par}) as equality, then
$$
\ba{rcl}
\tau & = & {\alpha \over \alpha + \mu \xi^*_0}, \quad 1 -
\tau \; = \; {\mu \xi^*_0 \over \alpha + \mu \xi^*_0}.
\ea
$$
Using the upper bound $A \geq \vp(x_0)$, we get the
following estimate:
$$
\ba{rcl}
\vp_{\mu}^* - \vp^* & \leq & \left[ {\mu \xi^*_0
	(1-\alpha)\over \alpha + \mu \xi^*_0} + \alpha \right] (A
- \vp^*) \; = \; {\mu \xi^*_0 + \alpha^2 \over \mu \xi^*_0
	+ \alpha} (A - \vp^*) .
\ea
$$
Denoting $\beta = \mu \xi^*_0$, we can find the optimal
$\alpha_*$ from the equation
$$
\ba{rcl}
{2 \alpha_* \over \beta + \alpha_*^2} = {1 \over \beta +
	\alpha_*}.
\ea
$$
Thus, $\alpha_* + \beta = \sqrt{\beta+\beta^2}$. This
means that $\alpha_* = {\beta \over \beta+\sqrt{\beta
		+\beta^2}}$. Hence,
$$
\ba{rcl}
{\beta (1- \alpha_*)  \over \beta + \alpha_*} + \alpha_*
\; = \; {2 \beta \over \beta + \sqrt{\beta+\beta^2}} \;
\leq \; 2 \sqrt{\beta}.
\ea
$$
In other words, we get the following bound:
\beq\label{eq-PBound}
\ba{rcl}
\vp_{\mu}^* - \vp^* & \leq & 2 \sqrt{\mu \xi^*_0}\,  (A -
\vp^*).
\ea
\eeq

Therefore, if we have an approximate solution $\bar{x}$
to the regularized problem, which satisfies~\eqref{eq-AppF},
we can ensure the following bound for the original problem 
$$
\ba{rcl}
\vp(\bar{x}) - \vp^{*}
& \overset{\eqref{eq-PBound}}{\leq} &
\vp(\bar{x}) - \vp_{\mu}^{*} + 2\sqrt{\mu \xi_0^*} (A - \vp^{*})
\;\; \overset{\eqref{eq-FUp}}{\leq} \;\;
\vp_{\mu}(\bar{x}) - \vp_{\mu}^{*} + 2\sqrt{\mu \xi_0^*} (A - \vp^{*}) \\
\\
& \overset{\eqref{eq-AppF}}{\leq} &
\delta( \vp(x_0) - \vp^{*}) + 2\sqrt{\mu \xi_0^*} (A - \vp^{*}),
\ea
$$
and the regularization parameter should be of the following order:
$$
\boxed{
\ba{rcl}
\mu & \approx & \frac{\delta^2}{\xi_0^{*}}
\ea
}
$$

Now, let us discuss a possible choice for the regularization functions.
Our goal is to have a uniformly convex smooth part of the objective. Thus, the following 
candidate for the regularizer is the most natural:
$$
\ba{rcl}
d_i(x) & := & f_i(x) + \frac{c_i}{p + 1}\|x - x_0\|^{p + 1},
\ea
$$
for a certain $c_i > 0$. This function is uniformly convex of degree $p + 1$
with parameter $\sigma_{p + 1}(d_i) = \frac{c_i}{2^{p - 1}}$
(see e.g. Lemma~2.5 in~\cite{doikov2021minimizing}).
Moreover, its $p$th 
derivative is Lipschitz continuous with constant 
$L_p(d_i) = L_p(f_i) + c_i  \cdot p! $ 
(see Theorem 7.1 in~\cite{rodomanov2020smoothness}).

Hence, applying method~\eqref{met-MFull}
to the regularized objective, we obtain the linear rate
$$
\ba{rcl}
\vp_{\mu}(x_k) - \vp_{\mu}^{*}
& \overset{\eqref{eq-RateF}}{\leq} &
(1 - \beta)^k (\vp_{\mu}(x_0) - \vp_{\mu}^{*}),
\ea
$$
where the condition number is equal to
\beq \label{CondNumber1}
\ba{rcl}
\beta & = & {\hat \beta}_p((1 - \mu)f + \mu d) 
\;\; = \;\; 
\min\limits_{1 \leq i \leq m} 
\frac{(p! \gamma_p((1 - \mu)f_i + \mu d_i))^{\frac{1}{p} } 
}{  (1 + p)^{\frac{1}{p}} + (p!\gamma_p((1 
- \mu)f_i + \mu d_i))^{\frac{1}{p}}  },
\ea
\eeq
with $\gamma_p((1 - \mu)f_i + \mu d_i) 
= \frac{\sigma_{p + 1}((1 - \mu)f_i + \mu d_i)}{L_p((1 - \mu)f_i + \mu d_i)}
= \frac{\mu c_i}{2^{p - 1}( L_p(f_i) + \mu c_i p! )}$.
We see that it is natural to set $c_i := L_p(f_i)$. In this case, we have
$$
\ba{rcl}
\beta & = & 
\Bigl[
1 + \bigl( 
(1 + p) 2^{p - 1}\bigl( \frac{1}{\mu p!} + 1 \bigr)
\bigr)^{\frac{1}{p}}
\Bigr]^{-1}.
\ea
$$
Thus, parameter $\mu$
plays a crucial role in the complexity of 
regularized problem~\eqref{prob-CompM}.

\section{Subhomogeneous Functions}
\label{sc-Subhomo}
\SetEQ

In this section, we consider a finer problem class
by adding some additional assumption on  
the outer component of the fully composite objective.
We show that for such problems, it is possible
to prove the global convergence rates
for the methods in a general convex case, when the smooth part is not necessary uniformly convex.
At the same time, we demonstrate that 
our methods can be accelerated.

A closed convex function $f: \dom f \to \R$ is called
\textit{subhomogeneous} 
if for any $x \in \dom f$ and $\gamma \geq 1$
such that $\gamma x \in \dom f$,
we have
\beq \label{DefSubhomo}
\ba{rcl}
f(\gamma x) & \leq & \gamma f(x).
\ea
\eeq
\BT 
Closed and convex function $f(\cdot)$ is subhomogeneous
if and only if it satisfies one of the following three conditions:
\begin{align}
\la g_x, x \ra  & \;\; \leq \;\;  f(x), 
\quad
\qquad x \in \dom f, \quad g_x \in \partial f(x), \label{Subhomo1} \\[3pt]
\la g_y, x \ra & \;\; \leq \;\; f(x), 
\quad
\qquad x, y \in \dom f, \quad g_y \in \partial f(y), \label{Subhomo2} \\[3pt]
f(x + t y) & \;\; \leq \;\; f(x) + t f(y),
\qquad x, y, x + ty \in \dom f, \quad t \geq 0. \label{Subhomo3}
\end{align}
\ET
\proof
Assume that~\eqref{DefSubhomo} is true. Then
$$
\ba{rcl}
\gamma f(x)
& \geq & f(\gamma x)
\;\; \geq \;\; f(x) + \la g_x, (\gamma - 1) x \ra,
\ea
$$
and this is~\eqref{Subhomo1}.

Assume~\eqref{Subhomo1} is true. Since $f$ is convex,
for any $x, y \in \dom f$ and $g_y \in \partial f(y)$, we have
$$
\ba{rcl}
f(x) - \la g_y, x \ra 
& \geq & f(y) - \la g_y, y \ra 
\;\; \overset{\eqref{Subhomo1}}{\geq} \;\; 0.
\ea
$$
This is relation~\eqref{Subhomo2}.

Finally, assume that~\eqref{Subhomo2} is true.
For any $y \in \dom f$ denote by $g(y)$ a particular
subgradient in $\partial f(y)$. Then, for any $x \in \dom f$
and $\gamma > 1$ such that $\gamma x \in \dom f$, 
we have
$$
\ba{rcl}
f(\gamma x)
& = & f(x) + \int\limits_{0}^{\gamma - 1}
\la g(x + \tau x), x \ra d\tau
\;\; \overset{\eqref{Subhomo2}}{\leq} \;\;
f(x) + \int\limits_{0}^{\gamma - 1} f(x) d\tau
\;\; = \;\; f(x).
\ea
$$
And this is~\eqref{DefSubhomo}.
Thus, conditions~\eqref{DefSubhomo}, \eqref{Subhomo1}
and \eqref{Subhomo2} are equivalent.

In order to justify equivalence with~\eqref{Subhomo3},
note that it can be rewritten as
$$
\ba{rcl}
\frac{1}{t}[ f(x + ty) - f(x) ] & \leq & f(y). 
\ea
$$
Therefore,
$$
\ba{rcl}
\max\limits_{g \in \partial f(x)} \la g, y \ra
& = & 
\lim\limits_{t \to +0}
\frac{1}{t}[ f(x + ty) - f(x) ] \;\; 
\leq \;\; f(y),
\ea
$$
and this is~\eqref{Subhomo2}.
On the other hand, if~\eqref{Subhomo2} is true, then
$$
\ba{rcl}
f(x + t y) - f(x) & = & \int\limits_{0}^{t} 
\la g(x + \tau y), y \ra d\tau
\;\; \overset{\eqref{Subhomo2}}{\leq} \;\;
\int\limits_{0}^t f(y) d\tau \;\; = \;\; t f(y).
\ea
$$
And this is~\eqref{Subhomo3}.
\qed

\BE
Clearly, function
$
f(x) \; = \; \max\limits_{i = 1}^n x^{(i)}
$
is subhomogeneous.
\EE

\BE
Consider the following function:
$$
\ba{rcl}
f(x) & = & \ln \biggl( \,
\sum\limits_{i = 1}^n e^{x^{(i)}}
\, \biggr)
\;\; = \;\;
\max\limits_{u \in \Delta_n}
\{ \la u, x \ra - \eta(u) \},
\ea
$$
where $\Delta_n \in \R^n_{+}$
is the standard simplex,
and $\eta(u) = \sum_{i = 1}^n u^{(i)} \ln u^{(i)}$
is the negative entropy.
Denote by $u(x)$ the unique optimal solution to this problem.
Then $\nabla f(x) = u(x)$, and we conclude that
$$
\ba{rcl}
\la \nabla f(x), x \ra & = &
\la u(x), x \ra \;\; \leq \;\;
\la u(x), x \ra - \eta(u(x)) \;\; = \;\; f(x),
\ea
$$
since $\eta(u) \leq 0$ for all $u \in \Delta_n$.
Thus, function $f(\cdot)$ is subhomogeneous since
condition~\eqref{Subhomo1} is satisfied.
\qed
\EE

\BE
Let function $f$ be subhomogeneous. Then, for a linear operator $A$,
function 
$$
\ba{rcl}
\bar{f}(x) & = & f(Ax)
\ea
$$ is subhomogeneous too.
Further, we can handle any affine transformation $Ax + b$, by
incorporating into our problem an auxiliary variable $\tau \in \R$:
$$
\ba{rcl}
\bar{f}(x, \tau) & = & f(Ax + \tau b),
\ea
$$
with additional normalizing constraint $\tau = 1$. \qed
\EE

\BR Let $F : \R^m \to \R \cup \{ +\infty \}$ be a closed convex function.
Assume that $F$ is monotone on its domain. Consequently, with any
point $x \in \dom F$, we have
$$
\ba{rcl}
x - \R^m_{+} & \subseteq & \dom F.
\ea
$$
Hence, if the domain contains a vector with strictly positive entries,
we have 
\beq \label{ZeroInDom}
\ba{rcl}
0 & \in & \inter( \dom F).
\ea
\eeq

Assume that for all $x, y \in \dom F$ and any $t \geq 0$, it holds
\beq \label{SubHomF}
\ba{rcl}
F(x + ty) & \leq F(x) + t F(y).
\ea
\eeq
\ER
Combining~\eqref{ZeroInDom} and~\eqref{SubHomF},
we conclude $\dom F = \R^m$. \qed

Now, let us introduce
our additional assumption.
\BAS
For any $x \in \dom \vp$,
the function $F(x, u)$ 
is subhomogeneous in $u \in \R^m$.
Thus, for any $x \in \dom \vp$, it holds
\beq \label{FSubHomo}
\ba{rcl}
F(x, u + t v) & \leq & F(x, u) + t F(x, v),
\quad
u, v \in \R^m, 
\quad t \geq 0.
\ea
\eeq
\EAS

We are ready to analyze convergence of the
Full-step $p$th order Basic Method~\eqref{met-MFull}
in a general convex case, when uniform convexity is absent. Let us use the following notation:
$$
\ba{rcl}
F(L_p(f)) & \Def & \sup\limits_{x \in \dom \varphi} F(x, L_p(f)).
\ea
$$

\BT \label{th-FullConvex}
Let the initial level set be bounded
\beq \label{LevelSet}
\ba{rcl}
D_0 & \Def &
\sup\limits_{x} 
\Bigl\{
\|x - x^{*}\| \; : \;
\vp(x) \leq \vp(x_0)
\Bigr\}
\;\; < \;\; +\infty.
\ea
\eeq
Then, for the iterations $\{x_k\}_{k \geq 1}$
of the method~\eqref{met-MFull}, we have
\beq \label{MFullConvex}
\ba{rcl}
\vp(x_k) - \vp^{*} & \leq & \frac{ (p + 1)^{p + 1}
	F(L_p(f)) D_0^{p + 1}}{p!} \cdot k^{-p}.
\ea
\eeq
\ET
\proof
By the definition of the method step, we have
\beq \label{OneStepGuarantee}
\ba{rcl}
\vp(x_{k + 1}) & = & F(x_{k + 1}, f(x_{k + 1})) \\
\\
& \overset{\eqref{eq-BoundF}}{\leq} &
F\Bigl(x_{k + 1}, 
\Omega_p(f, x_k; x_{k + 1}) + \frac{pL_p(f)}{(p + 1)!}\|x_{k + 1} - x_k\|^{p + 1} \Bigr) \\
\\
& \leq &
F\Bigl(y, \Omega_p(f, x_k; y) + \frac{pL_p(f)}{(p + 1)!}\|y - x_k\|^{p + 1} \Bigr) \\
\\
& \overset{\eqref{eq-BoundF}}{\leq} &
F\Bigl(y, f(y) + \frac{L_p(f)}{p!}\|y - x_k\|^{p + 1} \Bigr) \\
\\
& \overset{\eqref{FSubHomo}}{\leq} &
\varphi(y) + \frac{F(L_p(f))}{p!}\|y - x_k\|^{p + 1},
\ea
\eeq
for all $y \in \dom \vp$. Substituting $y := x_k$, we conclude
$$
\ba{rcl}
\vp(x_{k + 1}) & \leq & \vp(x_k).
\ea
$$
Hence, the method is monotone.
Let us take a convex combination 
$y := \frac{a_{k + 1}}{A_{k + 1}}x^{*} + \frac{A_k}{A_{k + 1}} x_k$,
where 
$A_k := k \cdot (k + 1) \cdot \ldots \cdot (k + p)$, and
$a_{k + 1} := A_{k + 1} - A_k = \frac{p + 1}{k + p + 1} A_{k + 1} $.
Therefore, we obtain by convexity
\beq \label{GenOneStep}
\ba{rcl}
\vp(x_{k + 1}) & \leq &
\frac{a_{k + 1}}{A_{k + 1}} \vp^{*} 
+ \frac{A_k}{A_{k + 1}} \vp(x_k)
+ \bigl( \frac{a_{k + 1}}{A_{k + 1}} \bigr)^{p + 1}
\cdot \frac{F(L_p(f)) \|x^{*} - x_k\|^{p + 1}}{p!} \\
\\
& \leq &
\frac{a_{k + 1}}{A_{k + 1}} \vp^{*}
+ \frac{A_k}{A_{k + 1}} \vp(x_k) 
+ \bigl(  \frac{p + 1}{k + p + 1}  \bigr)^{p + 1}
\cdot \frac{F(L_p(f)) D_0^{p + 1}}{p!}.
\ea
\eeq
Multiplying both sides by $A_{k + 1}$, we get
$$
\ba{rcl}
A_{k + 1}( \vp(x_{k + 1}) - \vp^{*} )
& \leq &
A_k ( \vp(x_k) - \vp^{*} )
+ A_{k + 1} \bigl( \frac{p + 1}{k + p + 1} \bigr)^{p + 1}
\cdot \frac{F(L_p(f)) D_0^{p + 1}}{p!} \\
\\
& \leq &
A_k (\vp (x_k) - \vp^{*})
+ \frac{(p + 1)^{p + 1} F(L_p(f)) D_0^{p + 1}}{p!}.
\ea
$$
Summing up the last inequality for different iterations, 
we finally obtain, for $k \geq 1$:
$$
\ba{rcl}
\vp(x_k) - \vp^{*}
& \leq & \frac{1}{A_k} \cdot 
\frac{k (p + 1)^{p + 1} F(L_p(f)) D_0^{p + 1} }{p!}
\;\; \leq \;\;
\frac{(p + 1)^{p + 1} F(L_p(f)) D_0^{p + 1} }{p!}
\cdot k^{-p}. \QF
\ea
$$

\section{Fully Composite Gradient Methods}
\label{sc-Gradient}
\SetEQ

Let us consider a more efficient version
of the Fully Composite Methods for the particular case $p = 1$
(first-order algorithms).
We start with the basic scheme.

\beq\label{met-GM}
\ba{|c|}
\hline\\
\quad \mbox{\bf Basic Gradient Method} \quad\\
\\
\hline\\
\ba{l}
\mbox{{\bf Choose} $x_0 \in \dom \vp$ 
	and $M = \alpha F(L_1(f))$, $\alpha \geq 1$.} \\
\\
\mbox{\bf For $k \geq 0$ iterate:}\\[10pt]
x_{k+1} \; = \; \argmin\limits_{y \in \dom \vp}
\Bigl\{
F\bigl(y, f(x_k) +  \la \nabla f(x_k),y - x_k \ra  \bigr)
+ \frac{M}{2}\|y - x_k\|^2
\Bigr\}.
\\
\\
\ea\\
\hline
\ea
\eeq

Contrary to the scheme~\eqref{met-MFull},
the regularization term in the method~\eqref{met-GM}
is \textit{outside} of the composite part.
Therefore, an implementation of each step can be much simpler.

For one iteration of the method,
for all $y \in \dom \vp$, we have
$$
\ba{rcl}
\varphi(x_{k + 1})
& \overset{\eqref{eq-BoundF}}{\leq} &
F\Bigl(x_{k + 1}, 
f(x_k) + \la \nabla f(x_k), x_{k + 1} - x_k \ra + \frac{L_1(f)}{2}\|x_{k + 1} - x_k\|^2
\Bigr) \\
\\
& \overset{\eqref{FSubHomo}}{\leq} &
F\Bigl( x_{k + 1},
f(x_k) + \la \nabla f(x_k), x_{k + 1} - x_k \ra \Bigr) 
+ \frac{\alpha F(L_1(f))}{2}\|x_{k + 1} - x_k\|^2 \\
\\
& \overset{\eqref{met-GM}}{\leq} &
F\Bigl( y, f(x_k) + \la \nabla f(x_k), y - x_k \ra \Bigr)
+ \frac{\alpha F(L_p(f)) \|y - x_k\|^2}{2} \\
\\
& \leq &
\varphi(y) + \frac{\alpha F(L_p(f)) \|y - x_k\|^2}{2},
\ea
$$
where we used component-wise
convexity of $f$ and monotonicity of $F$
in the last inequality.
By the same arguments as in the proof of Theorem~\ref{th-FullConvex},
we get the following result.

\BT Let the initial level set be bounded~\eqref{LevelSet}.
Then, for the sequence $\{ x_k \}_{k \geq 1}$ generated by
the method~\eqref{met-GM}, it holds
$$
\ba{rcl}
\vp(x_k) - \vp^{*} & \leq & 
\frac{4 \alpha F(L_1(f)) D_0^2}{k}.
\ea
$$
\ET

For the problems with bounded domain, we can propose the following 
alternative scheme,
which is a generalization of the classical 
Frank-Wolfe algorithm~\cite{frank1956algorithm,nesterov2018complexity}.
In this method, we do not use an explicit regularizer.
Thus, the cost of each step is usually even cheaper than in the 
Gradient Method~\eqref{met-GM}.

\beq\label{met-CGM}
\ba{|c|}
\hline\\
\quad \mbox{\bf Contracting Conditional Gradient Method} \quad\\
\\
\hline\\
\ba{l}
\mbox{{\bf Choose} $x_0 \in \dom \vp$ and $\{ \gamma_k \}_{k \geq 0}$.} \\
\\
\mbox{\bf For $k \geq 0$ iterate:}\\[10pt]
x_{k+1} \; = \; \argmin\limits_{y}
\Bigl\{
F\bigl(y, f(x_k) +  \la \nabla f(x_k),y - x_k \ra  \bigr) \\
\qquad \qquad \qquad \;\;\;\;\;
\; : \; x_k + \frac{1}{\gamma_k}(y - x_k) \in \dom \vp
\Bigr\}.
\\
\\
\ea\\
\hline
\ea
\eeq

\BT Let $\dom \vp$ be a bounded convex set. Denote its diameter by
\beq \label{DiamDef}
\ba{rcl}
\D & \Def & \sup\limits_{x, y \in \dom \vp} \|x - y\|
\;\; < \;\; +\infty.
\ea
\eeq
Set $\gamma_k := \frac{2}{k + 2}$.
Then, for the iterations $\{ x_k \}_{k \geq 1}$
of the method~\eqref{met-CGM}, it holds
\beq \label{CGM_Rate}
\ba{rcl}
\vp(x_k) - \vp^{*} & \leq & 
\frac{4 F(L_1(f)) \D^2}{k}.
\ea
\eeq
\ET
\proof
Let us denote the point 
$v_{k + 1} \Def x_k + \frac{1}{\gamma_k}(x_{k + 1} - x_k)
\overset{\eqref{met-CGM}}{\in} \dom \vp$. Hence,
\beq \label{XkDist}
\ba{rcl}
\|x_{k + 1} - x_k\|
& = & \gamma_k \|v_{k + 1} - x_k\|
\;\; \leq \;\; \gamma_k \D.
\ea
\eeq
Now, considering one iteration of the method, we obtain
$$
\ba{rcl}
\vp(x_{k + 1}) 
& \overset{\eqref{eq-BoundF}}{\leq} &
F\Bigl(
x_{k + 1}, f(x_k) + \la \nabla f(x_k), x_{k + 1} - x_k \ra 
+ \frac{L_1(f)}{2}\|x_{k + 1} - x_k\|^2 \Bigr) \\
\\
& \overset{\eqref{FSubHomo}}{\leq} &
F\Bigl( x_{k + 1}, f(x_k) + \la \nabla f(x_k), x_{k + 1} - x_k\ra \Bigr)
+ \frac{F(L_1(f))}{2}\|x_{k + 1} - x_k\|^2 \\
\\
& \overset{\eqref{XkDist}}{\leq} &
F\Bigl( x_{k + 1}, f(x_k) + \la \nabla f(x_k), x_{k + 1} - x_k\ra \Bigr)
+ \frac{\gamma_k^2 F(L_1(f)) \D^2}{2} \\
\\
& \overset{\eqref{met-CGM}}{\leq} &
F\Bigl(y, f(x_k) + \la \nabla f(x_k), y - x_k \ra \Bigr)
+ \frac{\gamma_k^2 F(L_1(f)) \D^2}{2} \\
\\
& \leq & \vp(y) + \frac{\gamma_k^2 F(L_1(f)) \D^2}{2},
\ea
$$
for every $y \in \dom \vp$.
Substituting $y := \gamma_k x^{*} + (1 - \gamma_k) x_k$,
and using convexity of $\vp$, we obtain
inequality very similar to~\eqref{GenOneStep}, 
from the proof of Theorem~\ref{th-FullConvex} for $p = 1$.
Hence, by the same arguments, we establish the rate~\eqref{CGM_Rate}.
\qed

Finally, we can present an accelerated method
(see~\cite{nesterov1983method,nesterov2018lectures,daspremont2021acceleration}).
\beq\label{met-FGM}
\ba{|c|}
\hline\\
\quad \mbox{\bf Fully Composite Fast Gradient Method} \quad\\
\\
\hline\\
\ba{l}
\mbox{{\bf Choose} $x_0 \in \dom \vp$ and $M = \alpha F(L_1(f))$,
	$\alpha \geq 1$. Set $v_0 = x_0$, $A_0 = 0$.} \\
\\
\mbox{\bf For $k \geq 0$ iterate:}\\[10pt]
\mbox{1. Find $a_{k + 1}$ from the equation $A_k + a_{k + 1} = a_{k + 1}^2$.  
	Set $A_{k + 1} = A_k + a_{k + 1}$.} \\[10pt]
\mbox{2. $y_{k} = \frac{a_{k + 1}v_k + A_k x_k}{A_{k + 1}}$.} \\[10pt]
\mbox{3. $x_{k+1} \; = \; \argmin\limits_{y}
\Bigl\{
F\bigl(y, f(y_{k}) +  \la \nabla f(y_{k}),y - y_{k} \ra  \bigr)
+ \frac{M}{2} \|y - y_{k}\|^2
\Bigr\}.$} \\[10pt]
\mbox{4. $v_{k + 1} \; = \; x_{k + 1} + \frac{A_k}{a_{k + 1}}(x_{k + 1} - x_k)$.} \\
\\
\ea\\
\hline
\ea
\eeq

\BT
For the sequence $\{x_k \}_{k \geq 1}$, generated by the method~\eqref{met-FGM}, it holds
\beq \label{FGMInduct}
\ba{rcl}
A_k \vp (x_k) + \frac{M}{2}\|x - v_k\|^2
& \leq &
A_k \vp (x) + \frac{M}{2}\|x - x_0\|^2, \qquad x \in \dom \vp.
\ea
\eeq
Consequently, in the case $x = x^*$, we get the following convergence guarantees:
\beq \label{FGMRate}
\ba{rcl}
\vp(x_k) - \vp^{*} & \leq & \frac{M \|x^{*} - x_0\|^2}{ 2A_k }
\;\; \leq \;\; \frac{2M \|x^{*} - x_0\|^2}{k^2}.
\ea
\eeq
\ET
\proof
Let us establish~\eqref{FGMInduct} by induction.
Assume that it holds for the current iterate, and consider the next step.
We fix an arbitrary $x \in \dom \vp$ and denote 
the convex combination $y := \frac{a_{k + 1} x + A_k x_k}{A_{k + 1}}$. Then
\beq \label{FGM_proof1}
\ba{cl}
& \frac{M}{2}\|x - x_0\|^2
+ A_{k + 1} \vp(x)
\;\; = \;\; \frac{M}{2}\|x - x_0\|^2 + A_k \vp(x) + a_{k + 1} \vp(x) \\
\\
& \overset{\eqref{FGMInduct}}{\geq} \;
\frac{M}{2} \|x - v_k\|^2 + A_k \vp(x_k) + a_{k + 1} \vp(x) \\
\\
& \;\, \geq \;\;\,
\frac{M}{2}\|x - v_k\|^2 + A_{k+1} \vp(y)
\;\; = \;\; 
A_{k + 1} \Bigl( \frac{M}{2}\|y - y_{k}\|^2 + \vp(y)  \Bigr) \\
\\
& \;\, \geq \;\;\,
A_{k + 1} \Bigl( \frac{M}{2}\|y - y_{k} \|^2
+ F\bigl(y, f(y_{k}) + \la \nabla f(y_{k}), y - y_{k} \ra \bigr)
\Bigr),
\ea
\eeq
where in the last inequality we used convexity of components of $f$
and monotonicity of $F$. 

The function in the right hand side of~\eqref{FGM_proof1}
is \textit{strongly convex} in $y$. Hence, we obtain
$$
\ba{cl}
& \frac{M}{2}\|y - y_{k}\|^2 + 
F\bigl(y, f(y_{k}) + \la \nabla f(y_{k}), y - y_{k} \ra \bigr) \\
\\
& \overset{\eqref{met-FGM}}{\geq} \; 
\frac{M}{2}\|y - x_{k + 1}\|^2 
+ \frac{M }{2}\|x_{k + 1} - y_k\|^2 
+ F\bigl(x_{k + 1}, f(y_{k}) + \la \nabla f(y_{k}), x_{k + 1} - y_{k} \ra \bigr) \\
\\
& \overset{\eqref{FSubHomo}}{\geq} \;
\frac{M}{2}\|y - x_{k + 1}\|^2 
+
F\Bigl(x_{k + 1}, f(y_{k}) + \la \nabla f(y_{k}), 
x_{k + 1} - y_{k} \ra + \frac{L_1(f)}{2}\|x_{k + 1} - y_{k}\|^2   \Bigr) \\
\\
& \overset{\eqref{eq-BoundF}}{\geq} \;
\frac{M}{2}\|y - x_{k + 1}\|^2 + \vp(x_{k + 1})
\;\; = \;\;
\frac{M}{2A_{k + 1}} \|x - v_{k + 1}\|^2 + \vp(x_{k + 1}).
\ea
$$
Thus we establish~\eqref{FGMInduct} for all $k \geq 0$.

Note that
$\sqrt{A_{k + 1}} - \sqrt{A_k}
= 
\frac{a_{k + 1}}{\sqrt{A_{k + 1}} + \sqrt{A_k}}
\;\; \geq \;\; \frac{a_{k + 1}}{2 \sqrt{A_{k + 1}}}
\;\; = \;\; \frac{1}{2}$.
Hence, $A_{k} \geq \frac{k^2}{4}$, and this proves
the rate of convergence~\eqref{FGMRate}.
\qed

\section{Fully Composite Newton Methods}
\label{sc-Newton}
\SetEQ

In this section, we analyze different variants
of the second-order methods for the fully composite 
formulation. Let us start with cubic regularization of the Newton's method~\cite{nesterov2006cubic}.

\beq\label{met-CN}
\ba{|c|}
\hline\\
\quad \mbox{\bf Fully Composite Cubic Newton Method} \quad\\
\\
\hline\\
\ba{l}
\mbox{{\bf Choose} $x_0 \in \dom \vp$ 
	and $M = \alpha F(L_2(f))$, $\alpha \geq 1$.} \\
\\
\mbox{\bf For $k \geq 0$ iterate:}\\[10pt]
x_{k+1} \; = \; \argmin\limits_{y \in \dom \vp}
\Bigl\{
F\bigl(y, \Omega_2(f, x_k; y) \bigr)
+ \frac{M}{6}\|y - x_k\|^3
\Bigr\}.
\\
\\
\ea\\
\hline
\ea
\eeq
Iterations of this type are simpler than those in the general method~\eqref{met-MFull}
with $p = 2$ since the regularization term is \textit{outside} the composite part.
At the same time, for any $y \in \dom \vp$, we have the following guarantee:
\beq \label{CubicNewtonGuarantee}
\ba{rcl}
\vp(x_{k + 1}) 
& \overset{\eqref{eq-BoundF}}{\leq} &
F\bigl(x_{k + 1}, \Omega_2(f, x_k; x_{k + 1}) 
+ \frac{L_2(f)}{6}\|x_{k + 1} - x_k\|^3 \bigr) \\
\\
& \overset{\eqref{FSubHomo}}{\leq} &
F\bigl( x_{k + 1}, \Omega_2(f, x_k; x_{k + 1}) \bigr)
+ \frac{M}{6}\|x_{k + 1} - x_k\|^3 \\
\\
& \overset{\eqref{met-CN}}{\leq} &
F\bigl(y, \Omega_2(f, x_k; y) \bigr) + \frac{M}{6}\|y - x_k\|^3 \\
\\
& \overset{\eqref{eq-BoundF}}{\leq} &
F\bigl(y, f(y) + \frac{L_2(f)}{6}\|y - x_k\|^3 \bigr)
+ \frac{M}{6}\|y - x_k\|^3 \\
\\
& \overset{\eqref{FSubHomo}}{\leq} &
\vp(y) + \frac{(1 + \alpha) F(L_2(f))}{6}\|y - x_k\|^3.
\ea
\eeq
Hence, using the same reasoning as in Theorem~\ref{th-FullConvex}, we
prove the global convergence result.

\BT Let the initial level set be bounded~\eqref{LevelSet}. Then,
for the sequence $\{ x_k \}_{k \geq 1}$ generated by the method~\eqref{met-CN}, we have
$$
\ba{rcl}
\vp(x_k) - \vp^{*} & \leq & 
\frac{ 9 (1 + \alpha) F(L_2(f)) D_0^3}{2 k^2}.
\ea
$$
\ET

\bigskip

In papers~\cite{doikov2020convex,doikov2020affine},
we looked at another modification of the Newton's method,
based on the {\em contracting idea}.
Let us present the corresponding version of the algorithm
for the Fully Composite Formulation.
This method can be seen as a second-order counterpart
of the Conditional Gradient Method~\eqref{met-CGM}.

\beq\label{met-ContrNewton}
\ba{|c|}
\hline\\
\quad \mbox{\bf Fully Composite Contracting Newton Method} \quad\\
\\
\hline\\
\ba{l}
\mbox{{\bf Choose} $x_0 \in \dom \vp$ and $\{ \gamma_k \}_{k \geq 0}$.} \\
\\
\mbox{\bf For $k \geq 0$ iterate:}\\[10pt]
x_{k+1} \; = \; \argmin\limits_{y}
\Bigl\{
F\bigl(y, \Omega_2(f, x_k; y)  \bigr) \\
\qquad \qquad \qquad \;\;\;\;\;
\; : \; x_k + \frac{1}{\gamma_k}(y - x_k) \in \dom \vp
\Bigr\}.
\\
\\
\ea\\
\hline
\ea
\eeq

Repeating the previous reasoning, we obtain the following result.

\BT Let the size of $\dom \vp$ be bounded  
by diameter $\D$~\eqref{DiamDef}.
Define $\gamma_k := \frac{3}{k + 3}$. Then,
for the sequence $\{x_k\}_{k \geq 1}$, generated by the method~\eqref{met-ContrNewton},
we have
$$
\ba{rcl}
\vp(x_k) - \vp^{*} & \leq & \frac{9 F(L_2(f)) \D^3}{k^{2}}.
\ea
$$
\ET

\section{Fully Composite Contracting Proximal Scheme}
\label{sc-ContrProx}
\SetEQ

In this section, we develop an accelerated second-order method.

The cubic regularization of the Newton's method 
was accelerated in~\cite{nesterov2008accelerating},
using the \textit{Estimating Functions} technique.
It is based on accumulating the gradients at the new points of the 
optimization process into a global linear model.
However, for the fully composite problems, we can guarantee 
only the progress in terms of the \textit{objective function},
and the good properties of the gradients are not easily available.

Therefore, we use inexact Contracting Proximal-Point iterations 
(see \cite{doikov2020contracting,guler1992new,lin2015universal})
as the basis of our accelerated scheme.
For simplicity, we consider the case $p = 2$
(second-order methods). Generalization to arbitrary $p \geq 1$ 
is more or less straightforward (see also~\cite{doikov2020inexact}).

Let us choose a prox-function, suitable for our problem class:
$$
\ba{rcl}
d(x) & := & \frac{\alpha}{3}\|x - x_0\|^3,
\qquad x, x_0 \in \E,
\qquad
\alpha \;\; := \;\; F(L_2(f)).
\ea
$$
It is well known that this function is \textit{uniformly convex} of
degree $3$ (see e.g.~\cite{doikov2021minimizing}), so it holds:
\beq \label{UnifConv}
\ba{rcl}
\rho_d(x; y) & \Def &
d(y) - d(x) - \la \nabla d(x), y - x \ra
\;\; \geq \;\;
\frac{\alpha}{6}\|y - x\|^3,
\qquad x, y \in \E.
\ea
\eeq
We need the following facts on Bregman divergence.
They can be checked in a direct way.
\begin{enumerate}
	\item For a closed convex function $\psi: \E \to \R \cup \{+\infty\}$ and a fixed prox center $v \in \E$, denote $h(x) := \psi(x) + \rho_d(v; x)$. 
Then, for the optimum $T := \argmin_x h(x)$, it holds
	\beq \label{StrongConvex}
	\ba{rcl}
	h(x) & \geq & h(T) + \rho_d(T; x), \qquad x \in \E.
	\ea
	\eeq
	In other words, the function $h(\cdot)$ 
	is \textit{strongly convex with respect to} $d(\cdot)$.
	
	\item For any $\bar{v}, v \in \E$ and $x \in \E$, it holds
	\beq \label{BregmanManipulation}
	\ba{rcl}
	\rho_d(\bar{v}; x) & = & \rho_d(v; x)
	+ \rho_d(\bar{v}; v)
	+ \la \nabla d(v) - \nabla d(\bar{v}), x - v \ra.
	\ea
	\eeq
\end{enumerate}

We use this prox function in the following general method.

\beq\label{met-AC}
\ba{|c|}
\hline\\
\quad \mbox{\bf Fully Composite Contracting Proximal-Point Scheme} \quad\\
\\
\hline\\
\ba{l}
\mbox{{\bf Choose} $x_0 \in \dom \vp$ and $\delta > 0$.
	Set $v_0 = x_0$, $A_0 = 0$.} \\
\\
\mbox{\bf For $k \geq 0$ iterate:}\\[10pt]
\mbox{1. Choose $A_{k + 1} = \bigl( \frac{k + 1}{3} \bigr)^3$.
Set $\gamma_k = \frac{A_{k + 1} - A_k}{A_{k + 1}}$.} \\[10pt]
\mbox{2. Form the subproblem $h_{k + 1}(x) = A_{k + 1} \vp(\gamma_k x + (1 - \gamma_k)x_k) + \rho_d(v_k; x)$.} \\[10pt]
\mbox{3. Find a point $v_{k + 1}$ s.t. $h_{k + 1}(v_{k + 1}) - h_{k + 1}^{*} \leq \delta$ 
	by the basic method~\eqref{met-CN}.} \\[10pt]
\mbox{4. $x_{k + 1} = \gamma_k v_{k + 1} + (1 - \gamma_k) x_k$.} \\
\\
\ea\\
\hline
\ea
\eeq

Let us justify first its rate of convergence.
Then we discuss the efficiency of implementation of the Step 3.

\BT
For the iterations of method~\eqref{met-AC}, we have
\beq \label{MainGuarantee}
\ba{rcl}
A_{k} \varphi(x_{k})
+ \rho_d(v_k; x)
& \leq & A_{k} \varphi(x)
+ \rho_d(x_0; x) + C_k(x),
\qquad x \in \E,
\ea
\eeq
where 
$$
\ba{rcl}
C_k(x) & \Def & 
k \delta + 
c_1(x) \cdot \sum\limits_{i = 1}^k \|x - v_i\|
+
c_2 \cdot \sum\limits_{i = 1}^k \|x - v_i\|^2,
\ea
$$
with
$c_1(x) \Def
\alpha^{1/3} (6 \delta)^{2/3}  + 2 \alpha^{2/3} (6 \delta)^{1/3} \cdot \|x - x_0\|$ and
$c_2 \Def 2 \alpha^{2/3}(6\delta)^{1/3}$.
\ET
\proof
Let us denote by $v_{k + 1}^{*}$ the \textit{exact}
minimizer of $h_{k + 1}(\cdot)$:
$v_{k + 1}^{*} \Def \argmin\limits_{x} h_{k + 1}(x)$.
Then, by the presence in the objective 
of the uniformly convex term, we get
$$
\ba{rcl}
\delta & \geq & h_{k + 1}(v_{k + 1}) - h_{k + 1}(v_{k + 1}^{*})
\;\; \overset{\eqref{StrongConvex}}{\geq} \;\; \rho_d(v_{k + 1}^{*}; v_{k + 1}) \\
\\
& \overset{\eqref{UnifConv}}{\geq} &
\frac{\alpha}{6}\|v_{k + 1} - v_{k + 1}^{*}\|^3.
\ea
$$
Hence, we can bound the distance between the exact minimizer 
and the approximate point $v_{k + 1}$ obtained by the inner method,
as follows:
\beq \label{VkDiff}
\ba{rcl}
\| v_{k + 1} - v_{k + 1}^{*} \|
& \leq &
\bigl(\frac{6 \delta}{\alpha}\bigr)^{\frac{1}{3}}.
\ea
\eeq

Now, assume that~\eqref{MainGuarantee} holds by 
for the current $k \geq 0$, and consider the next step of the method. 
Let us denote $a_{k + 1} := A_{k + 1} - A_k$. Then, for arbitrary $x \in \E$, we have
\beq \label{Induct1}
\ba{rcl}
\rho_d(x_0; x) + A_{k + 1} \varphi(x)
& = & 
\rho_d(x_0; x) + A_k \varphi(x) + a_{k + 1} \varphi(x) \\
\\
& \overset{\eqref{MainGuarantee}}{\geq} &
\rho_d(v_k; x) 
+ A_k \varphi(x_k)
+ a_{k + 1} \varphi(x) 
- C_k(x) \\
\\
& \geq &
\rho_d(v_k; x) + A_{k + 1} \varphi( \gamma_k x + (1 - \gamma_k) x_k)
- C_k(x) \\
\\
& = &
h_{k + 1}(x) - C_k(x) \\
\\
& \overset{\eqref{StrongConvex}}{\geq} &
\rho_d(v_{k + 1}^{*}; x)
+
h_{k + 1}(v_{k + 1}^{*}) - C_k(x),
\ea
\eeq
where the last inequality follows from {\em strong convexity} of
$h_{k + 1}(\cdot)$
{\em with respect to} $d(\cdot)$.

We can continue using our inexact solution $v_{k + 1}$ to
the auxiliary problem:
$$
\ba{rcl}
\rho_d(v_{k + 1}^{*}; x) + h_{k + 1}(v_{k + 1}^{*})
& \overset{\text{Step 3}}{\geq} &
\rho_d(v_{k + 1}^{*}; x) + h_{k + 1}(v_{k + 1}) - \delta \\
\\
& \overset{\eqref{BregmanManipulation}}{=} &
\rho_d(v_{k + 1}; x) + h_{k + 1}(v_{k + 1}) - \delta
+ \rho_d(v_{k + 1}^{*}; v_{k + 1}) \\
\\
& & 
\qquad
+ \; \la \nabla d(v_{k + 1}) - \nabla d(v_{k + 1}^{*}), x - v_{k + 1} \ra \\
\\
& \geq &
\rho_d(v_{k + 1}; x) + h_{k + 1}(v_{k + 1}) - \delta \\
\\
& &
\qquad
- \; \| \nabla d(v_{k + 1}) - \nabla d(v_{k + 1}^{*}) \|_{*} \cdot \|x - v_{k + 1}\|. 
\ea
$$
Now, since
$\nabla^2 d(x) = \| x - x_0 \| B + {1 \over \| x - x_0 \|} B (x-x_0)(x-x_0)^*B \preceq 2 \alpha \|x - x_0\| B$ for all $ x \in \E$,
we can bound
the difference between the gradients, as follows:
$$
\ba{rcl}
\| \nabla d(v_{k + 1}) - \nabla d(v_{k + 1}^{*}) \|_{*}
& = &
\| \int\limits_{0}^1 \nabla^2 d(v_{k + 1}^{*} + \tau(v_{k + 1} - v_{k + 1}^{*}))
d \tau ) (v_{k + 1} - v_{k + 1}^{*}) \|_{*} \\[10pt]
& \leq &
2 \alpha \cdot \|v_{k + 1} - v_{k + 1}^{*}\| \cdot
\int\limits_{0}^1
\| v_{k + 1}^{*} + \tau(v_{k + 1} - v_{k + 1}^{*}) - x_0 \| d \tau \\[10pt]
& \leq &
2 \alpha \cdot \|v_{k + 1} - v_{k + 1}^{*}\| \cdot
\Bigl(
\|v_{k + 1} - x_0\| + \frac{1}{2}\|v_{k + 1} - v_{k + 1}^{*}\| 
\Bigr) \\[10pt]
& \leq &
2\alpha \cdot \|v_{k + 1} - v_{k + 1}^{*} \|
\cdot \|x - v_{k + 1} \| \\[10pt]
& & 
\qquad
\; + \;
2\alpha \cdot \|v_{k + 1} - v_{k + 1}^{*} \|
\cdot \|x - x_0\| \\[10pt]
& & 
\qquad
\; + \;
\alpha \cdot \|v_{k + 1} - v_{k + 1}^*\|^2 \\[10pt]
& \overset{\eqref{VkDiff}}{\leq} &
2 \alpha^{2/3} (6 \delta)^{1/3} \cdot \|x - v_{k + 1}\|
\; + \; 
2 \alpha^{2/3} (6 \delta)^{1/3} \cdot \|x - x_0\| \\[10pt]
& &
\qquad
\; + \; \alpha^{1/3} (6 \delta)^{2/3}.
\ea
$$
Therefore, combining all component together, we conclude that
$$
\ba{rcl}
\rho_d(x_0; x) + A_{k + 1} \varphi(x)
& \geq &
\rho_d(v_{k + 1}; x) + h_{k + 1}(v_{k + 1}) - C_k(x) - \delta \\[10pt]
& &
\quad 
\; - \; \|x - v_{k + 1} \| \cdot  \Bigl(  
	\alpha^{1/3} (6 \delta)^{2/3}  + 2 \alpha^{2/3} (6 \delta)^{1/3} \cdot \|x - x_0\|
	\Bigr) \\[10pt]
& &
\quad
\; - \; 
\|x - v_{k + 1}\|^2 \cdot
	\Bigl(
	2 \alpha^{2/3}(6\delta)^{1/3}
	\Bigr) \\[10pt]
& \equiv &
\rho_d(v_{k + 1}; x) + h_{k + 1}(v_{k + 1})
- C_{k + 1}(x) \\[10pt]
& \geq &
\rho_d(v_{k + 1}; x) + A_{k + 1} \varphi(x_{k + 1}) - C_{k + 1}(x).
\ea
$$
Thus,~\eqref{MainGuarantee} is justified for all $k \geq 0$.
\qed

Let us substitute into~\eqref{MainGuarantee} $x=x^{*}$, and fix some $K \geq 0$. 
For all $0 \leq k \leq K$, we get
\beq \label{ResBound}
\ba{cl}
& \frac{\alpha}{6}\|x^{*} - v_k\|^3
\;\; \overset{\eqref{UnifConv}}{\leq} \;\; \rho_d(v_k; x^{*}) + A_k(\varphi(x_k) - \varphi^{*}) \\
\\
&  \overset{\eqref{MainGuarantee}}{\leq} \;\; \rho_d(x_0; x^{*})
+ K\delta
+ c_1(x^{*}) \sum\limits_{i = 1}^k \|x^{*} - v_i\|
+ c_2 \sum\limits_{i = 1}^k \|x^{*} - v_i\|^2
\;\; \Def \;\; R_k,
\ea
\eeq
and we need to estimate the quantities $R_k$ from above.
Note that 
\beq \label{RkRec}
\ba{rcl}
R_{k + 1} & = & R_k + c_1(x^{*}) \|x^{*} - v_{k + 1}\|
+ c_2 \|x^{*} - v_{k + 1}\| \\
\\
& \overset{\eqref{ResBound}}{\leq} &
R_k + a R_{k + 1}^{1/3} + b R_{k + 1}^{2/3},
\ea
\eeq
where
\beq \label{abConstDefs}
\ba{rcl}
a & \Def & c_1(x^{*}) \cdot \bigl( \frac{6}{\alpha} \bigr)^{1/3}
\;\; = \;\; 
6 \delta^{2/3} + 12 \delta^{1/3} \rho_d(x_0; x^{*})^{1/3}, \\
\\
b & \Def & c_2 \cdot \bigl( \frac{6}{\alpha} \bigr)^{2/3}
\;\; = \;\; 12 \delta^{1/3}.
\ea
\eeq
Dividing~\eqref{RkRec} by $R_{k + 1}^{1/3}$ and using monotonicity
of the sequence, we obtain quadratic inequality with respect to $R_{k + 1}^{1/3}$,
that is
$$
\ba{rcl}
\bigl( R_{k + 1}^{1/3} \bigr)^2
- b R_{k + 1}^{1/3} - (R_k^{2/3} + a) & \leq & 0.
\ea
$$
It can be resolved as follows:
$$
\ba{rcl}
R_{k + 1}^{1/3} & \leq &
\frac{b + \sqrt{b^2 + 4(R_k^{2/3} + a)}}{2}
\;\; \leq \;\; 
b + \sqrt{R_k^{2/3} + a}
\;\; \leq \;\;
R_k^{1/3} + b + \sqrt{a}.
\ea
$$
Hence, telescoping the last inequality, we get
$$
\ba{rcl}
R_k & \leq & \Bigl(  R_0^{1/3} + (b + \sqrt{a}) k  \Bigr)^3
\;\; \leq \;\; 9\Bigl( R_0 + k^3 b^3 + k^3 a^{3/2} \Bigr).
\ea
$$
Substituting the actual values of the parameters, we 
come to the following conclusion.
\BC
For the iterations of the method~\eqref{met-AC}, 
for all $k \geq 1$, it holds
\beq \label{ACGuarantee1}
\ba{rcl}
A_k (\vp(x_k) - \vp^{*}) + \rho_d(v_k; x^{*})
& \leq & 
O\Bigl(
\rho_d(x_0; x^{*}) + k^3 \delta + k^3 \delta^{1/2} \rho_d(x_0; x^{*})^{1/2} 
\Bigr).
\ea
\eeq
Hence,
$$
\ba{rcl}
\vp(x_k) - \vp^{*}
& \leq & O\Bigl( \frac{  \rho_d(x_0; x^{*})}{k^3} + \delta + \delta^{1/2} \rho_d(x_0; x^{*})^{1/2} \Bigr).
\ea
$$
And to solve the initial problem with $\varepsilon$-accuracy, we need to pick up
\beq \label{ACDelta}
\ba{rcl}
\delta & \approx & \min\Bigl\{ \varepsilon, \; \frac{\varepsilon^2}{\rho_d(x_0; x^{*})} \Bigr\}. \QF
\ea
\eeq
\EC

Let us apply now the Basic Cubic Newton Method~\eqref{met-CN} for solving the
subproblems at Step 3.
Denote the contracted smooth part by
$$
\ba{rcl}
\bar{f}(x) & := & f(\gamma_k x + (1 - \gamma_k) x_k),
\ea
$$
and the new outer part by
$$
\ba{rcl}
\bar{F}(x, u) & := & A_{k + 1} F(\gamma_k x + (1 - \gamma_k) x_k, u) + \rho_d(v_k; x).
\ea
$$
Hence, the objective in the subproblem can be represented as follows:
\beq \label{ACSubproblem}
\ba{rcl}
h_{k + 1}(x) & \equiv & \bar{F}(x, \bar{f}(x))
\;\; \to \;\; \min\limits_{x \in \dom \vp}
\ea
\eeq
Then iterations of the method~\eqref{met-CN} applied to~\eqref{ACSubproblem} are:
$$
\ba{rcl}
z_{t + 1} & := & \argmin\limits_{x \in \dom \vp}
\Bigl\{ \bar{F}(x, \Omega_2(\bar{f}, z_t; x)) + \frac{\bar{F}(L_2(\bar{f}))}{6}\|x - z_t\|^3  \Bigr\},
\qquad t \geq 0,
\ea
$$
and let us start with $z_0 := v_k$.

Note that
$L_2(\bar f) = \gamma_k^3 L_2(f)$.
Consequently,
\beq \label{FBarBound}
\ba{rcl}
\bar{F}(L_2(\bar{f}))
& \leq & A_{k + 1} F(L_2(\bar{f}))
\;\; \overset{\eqref{DefSubhomo}}{\leq} \;\;
\gamma_k^3 A_{k + 1} F(L_2(f)) \\
\\
& = & \frac{a_{k + 1}^3}{A_{k + 1}^2}F(L_2(f))
\;\; = \;\; 
\frac{ ((k + 1)^3 - k^3)^3  }{3^3 (k + 1)^6}F(L_2(f))
\;\; \leq \;\; 
F(L_2(f)).
\ea
\eeq

The guarantee~\eqref{CubicNewtonGuarantee} of one step ensures that, 
for any $x \in \dom \vp$, it holds
\beq \label{ACOneStep}
\ba{rcl}
h_{k + 1}(z_{t + 1}) & \leq & h_{k + 1}(x) + \frac{\bar{F}(L_2(\bar{f}))}{3}\|x - z_t\|^3 \\
\\
& \overset{\eqref{FBarBound}}{\leq} & 
h_{k + 1}(x) + \frac{F(L_2(f))}{3}\|x - z_t\|^3.
\ea
\eeq
Let $x = \tau v_{k + 1}^{*} + (1 - \tau) z_t$,
with $\tau := \frac{1}{\sqrt{6}}$ and $v_{k + 1}^{*}$ being the minimizer of~\eqref{ACSubproblem}. Then
$$
\ba{rcl}
h_{k + 1}(z_{t + 1}) & \leq & 
\tau h_{k + 1}^{*} + (1 - \tau) h_{k + 1}(z_t) 
+ \frac{\tau^3 F(L_2(f))}{3}\|v^{*}_{k + 1} - z_t\|^3 \\
\\
& \overset{\eqref{UnifConv}}{\leq} &
\tau h_{k + 1}^{*} + (1 - \tau) h_{k + 1}(z_t)
+ 2 \tau^3 \bigl( h_{k + 1}(z_t) - h_{k + 1}^{*} \bigr).
\ea
$$
Therefore,
$$
\ba{rcl}
h_{k + 1}(z_{t + 1}) - h_{k + 1}^{*}
& \leq &
\Bigl(1 - \tau + 2 \tau^3\Bigr) \cdot \bigl( h_{k + 1}(z_t) - h_{k + 1}^{*} \bigr) \\
\\
& = & \Bigl(1 - \frac{2}{3 \sqrt{6}} \Bigr) \cdot \bigl( h_{k + 1}(z_t) - h_{k + 1}^{*} \bigr) 
\;\; \leq \;\; \frac{3}{4}\bigl( h_{k + 1}(z_t) - h_{k + 1}^{*} \bigr).
\ea
$$
We see that our subsolver has a fast \textit{linear} rate of convergence,
which does not depend on any condition number. 
Let us estimate the residual after one step of the method.
Substituting $x= x^{*}$ (the solution to the original problem) into~\eqref{ACOneStep},
we get
$$
\ba{rcl}
h_{k + 1}(z_1) - h_{k + 1}^{*}
& \overset{\eqref{ACOneStep}}{\leq} &
h_{k + 1}(x^{*}) - h_{k + 1}^{*} + \frac{F(L_2(f))}{3}\|x^{*} - v_k\|^3 \\
\\
& \leq &
a_{k + 1} \vp^{*} + A_k \vp(x_k) + \rho_d(v_k; x^{*}) - h_{k + 1}^{*}
+ \frac{F(L_2(f))}{3}\|x^{*} - v_k\|^3 \\
\\
& \overset{(*)}{\leq} &
A_k(\vp(x_k) - \vp^{*}) + \rho_d(v_k; x^{*}) + 2 \rho_d(v_k; x^{*}) \\
\\
& \overset{\eqref{ACGuarantee1},\eqref{ACDelta}}{\leq} &
O\Bigl( \rho_d(x_0; x^{*}) + k^3 \varepsilon \Bigr),
\ea
$$
where we used in $(*)$ the uniform convexity of the prox-function and 
the following bound:
$$
\ba{rcl}
h_{k + 1}(x) & \geq & 
A_{k + 1}F(\gamma_k x + (1 - \gamma_k) x_k, f(\gamma_k x + (1 - \gamma_k) x_k)) \\
\\
& \geq & \min\limits_{y \in \dom \vp } A_{k + 1} F(y, f(y))
\;\; = \;\; A_{k + 1} \vp^{*}.
\ea
$$

Combining these bounds together, we come to the following final conclusion.
\BC
For solving the initial problem with $\varepsilon$-accuracy:
$$
\ba{rcl}
\vp(x_K) - \vp^{*} & \leq & \varepsilon,
\ea
$$
we need to perform 
$K = O\Bigl( \Bigl[  \frac{F(L_2(f))\|x_0 - x^{*}\|^3}{\varepsilon} \Bigr]^{1/3} \Bigr)$
iterations of the proximal-point scheme~\eqref{met-AC}.
At each iteration, it requires no more than
$$
\ba{rcl}
N & = & 
O\Bigl(1 + \log\Bigl[ 
\frac{ F(L_2(f))\|x_0 - x^{*}\|^3 }{\varepsilon}
\Bigr]\Bigr)
\ea
$$
steps of the basic method~\eqref{met-CN}.
\EC

We see that the price to pay for the level of generality is an additional 
logarithmic term in the final complexity estimate.
It remains to be an open theoretical question: whether we can
develop a \textit{direct} accelerated high-order method 
for Fully Composite Formulation,
which does not need inexact proximal iterations.
It would also help in constructing
\textit{optimal} high-order 
methods~\cite{monteiro2013accelerated,gasnikov2019optimal,jiang2019optimal,bubeck2019near},
matching the existing lower complexity bounds~\cite{arjevani2019oracle,nesterov2018lectures}.

Another interesting research direction is the development of
\textit{universal}~\cite{grapiglia2017regularized,grapiglia2020tensor}
and \textit{randomized}~\cite{cartis2018global,hanzely2020stochastic}
variants of the Fully Composite Methods.


\end{document}